\documentclass[12pt]{article}
\usepackage{graphicx,epsfig}
\usepackage{color}
\usepackage{amsmath,amssymb,mathrsfs,amsthm}
\usepackage{amsfonts}
\usepackage{multicol}
 \usepackage{graphicx}
\usepackage{color}
\usepackage{pdfcolmk}
\usepackage[nice]{nicefrac}
\usepackage{amsthm}
\usepackage{latexsym}
\usepackage{dsfont}
\usepackage{setspace}
\usepackage[a4paper, tmargin=1.3in, bmargin=1in, lmargin=1in, rmargin=1in, headheight=13.4pt]{geometry}
\usepackage{natbib}
\usepackage[font=small,labelfont=bf,labelsep=space]{caption}
\makeatletter

\onehalfspace
\date{}

\def\@citex[#1]#2{\if@filesw\immediate\write\@auxout{\string\citation{#2}}\fi
  \def\@citea{}\@cite{\@for\@citeb:=#2\do
    {\@citea\def\@citea{,\linebreak[0]\hskip0pt plus .2em}%
      \@ifundefined{b@\@citeb}%
    {{\bf ?}\@warning{Citation `\@citeb' on page \thepage\space undefined}}%
      \hbox{\csname b@\@citeb\endcsname}}}{#1}}

\newtheorem{rule-def}[theorem]{Rule}
\numberwithin{equation}{section}
\makeatletter

\begin{document}
\begin{center}

\textbf{\normalsize \Large{Analytical Approach For Solving Population Balances: A Homotopy Perturbation Method}}
\vspace{3mm}

\textbf{Gurmeet Kaur$^{a}$, Randhir Singh$^b$, Mehakpreet Singh$^{c}$, Jitendra Kumar$^{a}$, Themis Matsoukas$^{d}$}\\
\vspace{2mm}

\textit{$^{a}$Department of Mathematics, Indian Institute of Technology Kharagpur, Kharagpur, India}\\
\textit{$^{b}$Department of Mathematics, Birla Institute of Technology  Mesra, Ranchi-835215, India}\\
\textit{$^{c}$BIOMATH, Department of Mathematical Modelling, Statistics and Bioinformatics,\\
Faculty of Bioscience Engineering, Ghent University, 9000 Ghent, Belgium} \\
\textit{$^d$Department of Chemical Engineering, Pennsylvania State University, 158 Fenske Laboratory,
University Park, PA 16802-4400, USA}\\

\end{center}

\begin{abstract}
In the present work, a new approach is proposed for finding the analytical solution of population balances. This approach is relying on idea of Homotopy Perturbation Method (HPM). The HPM solves both linear and nonlinear initial and boundary value problems without nonphysical restrictive assumptions such as linearization and discretization. It gives the solution in the form of series with easily computable solution components. The outcome of this study reveals that the proposed method can avoid numerical stability problems which often characterize in general numerical techniques related to this area. Several examples including Austin's kernel available in literature are examined to demonstrate the accuracy and applicability of the proposed method.
\end{abstract}


\section{Introduction}
Fragmentation, along with aggregation, represents one of the most basic particulate processes. It refers to the process by which a granule, colloidal cluster, polymer chain or other ``particle'' breaks into a distribution of fragments. This generic process is encountered in a wide range of natural systems and serves as a model of branching processes in general, in which a population is transformed by replacing a member by a distribution of new members. The physical mechanisms that give rise to breakup or fragmentation (the two terms are used interchangeably in this paper), are generally different between systems but a common mathematical description is possible under the rather general assumption that the rate at which a particle breaks up and the distribution of fragments that are produced depend only on the parent particle and are independent by the presence of other particles in the population. This model views breakup as a first order reaction in the concentration of parent particles, with a distribution of products that depend on the parent particle. Moreover, the physical mechanism that give rise to the aggregation process can be described under the general assumption that the rate at which the particles aggregate and the distribution that are depend on the particles which aggregate together to form a new particle.

The breakage as well as aggregation models will be described by two functions. For defining the breakage model, the rate of breakup of the parent particle, $a(m)$, where $m$ is the size of the particle, and the distribution of fragments, $k(m|n)$ that gives the number of fragments of size $m$ produced by a parent of size $n$ are required. If we take size to refer to the mass of the particle, the fragment distribution obeys the following normalization conditions:
\begin{equation}
    \int_0^n k(m|n) dm = \bar f,
\end{equation}
\begin{equation}
    \int_0^n m b(m|n) dm = n.
\end{equation}
The first condition gives the average number of fragments, $\bar f$, and is the same for all parent sizes $m$; the second condition expresses mass conservation between the parent particle and the fragments. A physically realistic fragmentation model requires $\bar f\geq 2$. The governing equation for the size distribution of population that undergoes the above fragmentation process is given by the fragmentation equation,
\begin{equation}\label{pbebrk}
    \frac{\partial c(m,t)}{\partial t} =\underbrace{-a(m) c(m,t)}_{\text{death of particle of size $m$}} + \underbrace{\int_m^\infty a(n) k(m|n) c(n,t) dn}_{\text{birth of particle of size $m$}},
\end{equation}
where $c(m,t)$ is the concentration of particles whose mass is in $(m,m+dm)$ (we normalize the mass $m$ by its mean value at time 0).

Moreover, the aggregation model can be defined with the help of aggregation kernel, $g(m,m')$, where $m$ and $m'$ are the sizes of the particles which defines the rate at which particles aggregate. The governing equation to track the change in particle size distribution due to aggregation process is defined by aggregation population balance equation,
 \begin{align}\label{sec1:eq1}
\nonumber \frac{\partial c(t,m)}{\partial t}&=\underbrace{\frac{1}{2}\int\limits_{0}^{m} g(m-n, n)c(t,m-n)c(t,n) dn}_{\text{birth of particle of size $m$ due to aggregation of particles of sizes $m-n$ and $n$}}\\& -\underbrace{\int\limits_{0}^{\infty}g(m, n)c(t,m)c(t,n)dn}_{\text{death of particles of sizes $m$ due to aggregation of particles of sizes $m$ and $n$}},~~t\in [0,T], ~~x\in(0,\infty)
\end{align}
with the initial condition
\begin{align}\label{sec1:eq2}
 c(0,m)=\xi_0(m).
\end{align}
Here $c(t,m)$ represents the concentration of particles of size $m$ at time $t$. Here, $g(m,n)$ is the aggregation  kernel, which describes the rate at which the particles of sizes $m$ and $n$ coagulate to form a particle of size $m+n$. It is non-negative and symmetric.

The equations \eqref{pbebrk} and \eqref{sec1:eq1} have been the subject of several investigations and the analytical solutions have been obtained for a small number of special cases. This work, largely theoretical, \citep{Simha:JAP41,Simha:JCP56,Tobolsky:JPS57,Ziff:M86,Ernst:JPMG93,Singh:PR96} seeks analytical solutions to the fragmentation equation and studies the scaling behavior of the size distribution at long times. A review of these developments and a collection of analytical and solutions is given by Ziff \citep{Ziff:JPMG91}.

The analytical solution for aggregation equation can be studied by many researchers for limited number of aggregation kernels. \cite{ranjbar2010numerical} derived analytical solution for constant kernel using Taylor polynomials and radial basis functions. Moreover, \cite{hammouch2010laplace} uses the Laplace-Variational iteration to obtain approximate series solution for equation \eqref{sec1:eq1}. For additive aggregation kernel,  \cite{fernandez2007exact} proposed analytical solution. However, exact solutions for constant, sum and product kernels are given in \cite{deaconu2000smoluchowski}.



The method of homotopy perturbation has been recently developed to obtain analytical solutions from differential and integral equations and has been successfully used to to solve a a wide range of dynamical problems \citep{Ganji:PL06,He:AMC03,Shahed:IJNSNS11}. The purpose of this paper is to formulate the homotopy perturbation method (HPM) in a form appropriate for the fragmentation and aggregation equations and to use it for obtaining analytical solutions. The paper is organized as follows. In section \ref{s:hpm} we introduce the general methodology of HPM. In section \ref{s:sol}, HPM will be adapt it to the fragmentation and aggregation equations. In next sections \ref{s:brkmodel} and \ref{s:aggmodel}, we apply the method to a number of known and new fragmentation and aggregation models. Finally in section \ref{s:discussion} we examine the transient and scaling solution for a  quaternary breakup and aggregation models and relate its behavior.

\section{The Homotopy Perturbation Method}
\label{s:hpm}
Recently, the homotopy perturbation method (HPM) has been studied by many researchers for solving linear and nonlinear problems \cite{he1999homotopy,he2000coupling,he2006homotopy,nazari2013modified,he2003homotopy,el2005application,he2006homotopy,ganji2006application}. The HPM yields a very rapid convergence of the solution series in most cases, usually only few iterations leading to very accurate solutions.
To illustrate the HPM, we consider the following differential equation
\begin{align}\label{sec2:eq4}
T(u) - g(r)= 0,~~~~ r\in\Omega,
\end{align}
with  boundary conditions
\begin{align}\label{sec2:eq5}
B(u,\partial u/\partial n)= 0,~~~~ r\in \partial\Omega,
\end{align}
where $T$  is a general differential operator and  $B$ is a boundary operator. Usually the operator $T$ can be decomposed into two parts, a linear operator $L$ and a nonlinear operator $N$, and expressed as
\begin{align}\label{sec2:eq6}
L(u)+N(u) - g (r)= 0.
\end{align}
NOTE: If the equation is linear then the non-linear part will be zero, i.e., $N(u)=0$.\\
According to HPM, we construct a homotopy  which satisfies
\begin{align}\label{sec2:eq7}
H[v(r,p)]=(1-p)[L(v(r,p))-L(u_0)]+p[T[v(r,q)] - g(r)]= 0,
\end{align}
where  $u_0$ is an initial guess to exact solution of \eqref{sec2:eq4}. When $p=0$ then $L(v(r,0)) = L(u_0)=0$,
and when $p=1$, then $T(v(r,1))- g(r) = 0. $ It is worth noting that as the embedding parameter $p$ increases monotonically from zero to unity, the changing process of $p$ from zero to unity is just that of $v(r,p)$ from $u_0(r)$ to $u(r)$. This is called deformation, and also $L(v)-L(u_0)$ and $T(v)- g(r)$ are called homotopic in topology.

According to the HPM, we can first view the embedding parameter $p$ as a small parameter, and construct the solution as a power series in $p$, as
\begin{align}\label{sec1:eq8}
   u=\sum_{k=0}^{\infty} p^{k}v_k=v_0+p v_1+p^2 v_2+ \cdots.
\end{align}
Substituting \eqref{sec1:eq8} in  \eqref{sec2:eq7} and then  letting $p=1$, we obtain the solution as
\begin{align}\label{sec1:eq9}
\displaystyle
 f=\lim_{p\rightarrow1} u=\sum_{k=0}^{\infty}v_k.
\end{align}
The series  \eqref{sec1:eq9} is a convergent for most of the cases  and the rate of convergence depends on the nature of the problem \cite{he2000coupling}.

The idea of the HPM for different types of models \citep{he2003homotopy,el2005application,he2006homotopy,ganji2006application} will be demonstrated in this section. Firstly, the problem linear integral equation will be solved using HPM. Secondly, non-linear integral equation will be taken into account.\\

\textbf{Model 1:} Consider the linear integral equation:
\begin{align}\label{model1}
c(t)=a(t)+\int_q^t K(u,t,c(u))du.
\end{align}
Here $K(u,t,c(u))$ should be linear function, i.e., $K(u,t,c(u))=\sum_{j=1}^{n}K_j(u,t)c_{i,j}(u)$. For solving above equation using HPM method, a homotopy will be constructed
\begin{align}\label{homotopy}
(1-p)(C(t)-c_0(t))+p(C(t)-a(t)-\int_q^t K(u,t)C(u)du)=0.
\end{align}
Suppose above equation has the solution of the following form:
\begin{align}\label{solser}
C(t)=C_0(t)+pC_1(t)+p^2 C_2(t)+\hdots,
\end{align}
where $C_i(t),~i=1,2,3,\hdots$ are the functions which will be determined and $C_0(t)$ is the initial approximation to the solution which is taken to be $a(t)$. On substituting \eqref{solser} into \eqref{homotopy} and equating the like powers of $p$, we have
\begin{align*}
p^0 &:~ C_0(t)-c_0(t)=0,\\
p^1 &:~ C_1(t)+c_0(t)-a(t)-\int_q^t \sum_{j=1}^{n} K_j(u,t) C_0(u)du=0\\
p^2 &:~ C_2(t)-\int_q^t \sum_{j=1}^{n} K_j(u,t)C_0(u)du=0\\
\vdots
\end{align*}

The approximated solution of \eqref{model1} can be obtained by setting $p=1$,
\begin{align*}
c_i(t)=\lim_{p \to 1} C(t)  = \sum_{j=0}^{\infty} C_j(t).
\end{align*}

\textbf{Model 2:} Consider the nonlinear  integral equation
\begin{align}\label{sec1:eq4}
u(x)=g(x)+\int_{a}^{b} k(x,s) f(s,u(s)) ds,~~x\in\Omega.
\end{align}
To apply the HPM, we write \eqref{sec1:eq4} as
\begin{align}\label{sec1:eq5}
L(u)=u(x)-g(x)-\int_{a}^{b} k(x,s) f(s,u(s)) ds=0,
\end{align}
with solution $y(x)$ and we construct the homotopy $H(u, p)$ as
\begin{align}\label{sec1:eq6}
   H(u,0) = F(u),~~~~H(u,1)=L(u),
\end{align}
where $F(u)$ is a functional operator with known solution $v_0$. We choose a convex homotopy by
\begin{align}\label{conhomo}
   H(u,p) = (1-p)F(u)+pL(u)= 0,~~~p\in [0,1]
\end{align}
and continuously trace an implicitly defined curve from a starting point $H(v_0, 0)$ to
the solution function $H(f, 1)$. The embedding parameter $p$ increases monotonically from $0$ to $1$ as a trivial problem $F(u)=0$ is continuously deformed to the original problem $L(u)= 0$.

According to the HPM, we can first view the embedding parameter $p$ as a small parameter, and construct the solution as a power series in $p$, as
\begin{align}\label{useries}
   u=\sum_{k=0}^{\infty} p^{k}v_k.
\end{align}
 Substituting \eqref{useries} in  \eqref{conhomo} and then  letting $p=1$, we obtain the solution as
\begin{align}\label{ysol}
\displaystyle
 y=\lim_{p\rightarrow1} u=\sum_{k=0}^{\infty}v_k.
\end{align}
The series  \eqref{ysol} is a convergent for most of the cases  and the rate of convergence depends on the nature of the problem \cite{he2000coupling}.

Substitute \eqref{useries} into \eqref{sec1:eq5} and equate the terms with identical powers of $p$, we obtain
\begin{align}\label{sec3:eq5}
&u_0=g(x),\\
&u_{k+1}=\int_{a}^{b} k(x,s) H_n ds,~~~k=0,1,2\ldots
\end{align}

\section{Implementation of HPM}
\label{s:sol}
In this section, HPM will be implemented to the aggregation and fragmentation equations. It can be noticed clearly that the breakage equation is linear whereas aggregation equation is non-linear. Therefore, for solving breakage equation, HPM for linear equation will be implemented and for solving aggregation equation, HPM for non-linear equation will be implemented.

\subsection{Implementation of HPM for Breakage Equation}
To apply HPM to fragmentation we first express the population balance equation in the integral form,
\begin{equation}
\label{pbe:integral}
    \frac{\partial c(t,m)}{\partial t} =
       - a(m)c(t,m)
       + \int\limits_m^\infty k(m|n)a(n) c(t,n) dn,
\end{equation}
where $c_0(m) = c(t,m)$ is the size distribution at time zero. Next we introduce a new function $u=u(t,m; p)$ of time, size, $p$ and define the homotopy of Eq.\ (\ref{pbe:integral}) as follows:
\begin{multline}
\label{homo}
    (1-p) \left(\frac{\partial c(t,m)}{\partial t}-\frac{\partial c_0(t,m)}{\partial t}\right) +
    p\left\{\vphantom{\int\limits_m^\infty}
       \frac{\partial c(t,m)}{\partial t}
       +a(m)c(t,m) + \right.
       \\
       \left.
       -\int\limits_m^\infty k(m|n)a(n) c(t,n)dn
       \right\} = 0 .
\end{multline}
With $p=0$, Eq.\ (\ref{homo}) gives $u(t,m,0)=c(0,m)$; with $p=1$ it gives $u(t,m; 1) = c(t,m)$. Thus by continuously varying $p$ from 0 to 1 we obtain a continuous transformation of $u(t,m;p)$ from the initial distribution $c_0$ to the actual distribution $c(m,t)$ at time $t$.
Following HPM, we express $u(t,m;p)$ as a power series in $p$,
\begin{equation}
\label{series:p}
u(t,m;p) = \sum_{k=0}^\infty u_k(t,m) p^k,
\end{equation}
where the coefficients $u_k $ are function of time and size to be determined. According to the homotopy in Eq.\ (\ref{homo}) the solution to the fragmentation equation is obtained by setting $p=1$ in the above series:
\begin{equation}
\label{series}
    c(m,t) = \lim_{p\to 1} u(t,m;p) = \sum_{k=0}^\infty u_k(t,m).
\end{equation}
The series is convergent and the rate of convergence depends on the nature of Eq.\ (\ref{pbe:integral}). To obtain the coefficients $u_k$ we insert the series  \eqref{series:p} into (\ref{homo}):
\begin{multline}
    (1-p)\left(\frac{\partial}{\partial t} \sum_{k=0}^\infty p^k u_k- \frac{\partial c_0}{\partial t}\right)
    + p
    \left\{\frac{\partial}{\partial t} \sum_{k=0}^\infty p^k u_k+
        a(m)\left(\sum_{k=0}^\infty p^k u_k\right)
       \right. -
       \\
    \left.
       \int\limits_m^\infty k(m|n) a(n)
       \left(\sum_{k=0}^\infty p^k u_k\right) dn
    \right\}
    = 0  ,
\end{multline}
collect terms in powers of $p$ and set their coefficients to zero. We find:
\begin{align}
\label{p0}
   k=0 &  && u_0 = c(0,m) ,\\
\label{p1}
   k=1&  && \frac{\partial u_1}{\partial t} =
      -a(m) u_0(t,m) + \int_m^\infty k(m|n) a(n) u_0(t,n) dn, \\
   \vdots~~~~~& &&
   \nonumber\\
\label{pk}
   k~~~~&     && \frac{\partial u_k}{\partial t}=
   -a(m) u_{k-1}(t,m) + \int_m^\infty k(m|n) a(n) u_{k-1}(t,n) dn.
\end{align}
The first term, $u_0$, is equal to the initial distribution and each subsequent term involves the coefficient of previous order. Thus we may obtain these coefficients recursively. The $n$th order approximation of the solution is obtained by truncating the series past the $n$th term. Here, however, we are interested in cases for which the infinite series can be expressed in closed analytical form. First we reproduce known solutions to the fragmentation equation using HPM, and we also present new results for a case that has not been solved previously.

\subsection{Implementation of HPM for Aggregation Equation}
The homotopy of equation \eqref{sec1:eq1} is constructed which satisfies

\begin{align}
\nonumber H_p&=(1-p)\left(\frac{\partial u(t,m)}{\partial t}-\frac{\partial u_0(t,m)}{\partial t}\right)\\&+p\left(\frac{\partial u(t,m)}{\partial t}-\frac{1}{2}\int\limits_{0}^{m} a(m-n, n)u(t,m-n)u(t,n) dn+ \int\limits_{0}^{\infty}a(m, n)u(t,m)u(t,n)dn \right)=0
\end{align}

\begin{align}
\nonumber H_p&=\frac{\partial u(t,m)}{\partial t}-(1-p)\frac{\partial u_0(t,m)}{\partial t}\\&+p\left(-\frac{1}{2}\int\limits_{0}^{m} g(m-n, n)u(t,m-n)u(t,n) dn+ \int\limits_{0}^{\infty}a(m, n)u(t,m)u(t,n)dn \right)=0
\end{align}

Inserting $u=\sum_{k=0}^{\infty} p^{k}u_k$ in the above equation,
\begin{align}
\nonumber &\frac{\partial\sum_{k=0}^{\infty} p^{k}u_k(t,m)}{\partial t}-(1-p)\frac{\partial u_0(t,m)}{\partial t}-\frac{p}{2}\int\limits_{0}^{m} g(m-n, n)\sum_{k=0}^{\infty} p^{k}u_k(t,m-n)\sum_{k=0}^{\infty} p^{k}u_k(t,n) dn\\&+p\int\limits_{0}^{\infty}g(m, n)\sum_{k=0}^{\infty} p^{k}u_k(t,m)\sum_{k=0}^{\infty} p^{k}u_k(t,n)dn=0
\end{align}

Equate the terms with identical powers of $p$, we obtain the expression of HPM for aggregation equation:
\begin{align}
\nonumber &\frac{\partial u_1(t,m)}{\partial t}-\frac{\partial u_0(t,m)}{\partial t}-\frac{1}{2}\int\limits_{0}^{m} g(m-n, n) u_0(t,m-n)u_0(t,n) dn\\&+\int\limits_{0}^{\infty}g(m, n)u_0(t,m)u_0(t,n)dn=0\\
\nonumber &\frac{\partial u_2(t,m)}{\partial t}-\frac{1}{2}\int\limits_{0}^{m} g(m-n, n) [u_0(t,m-n)u_1(t,n)+u_1(t,m-n)u_0(t,n)] dn\\&+\int\limits_{0}^{\infty}g(m, n)[u_0(t,m)u_1(t,n)+u_1(t,m)u_0(t,n)]dn=0
\end{align}
\begin{align}
\nonumber &\frac{\partial u_3(t,m)}{\partial t}-\frac{1}{2}\int\limits_{0}^{m} g(m-n, n) [u_0(t,m-n)u_2(t,n)+u_2(t,m-n)u_0(t,n)+u_1(t,m-n)u_1(t,n)] dn\\&+\int\limits_{0}^{\infty}g(m, n)[u_2(t,m)u_0(t,n)+u_0(t,m)u_2(t,n)+u_1(t,m)u_1(t,n)]dn=0
\end{align}
\vdots
\begin{align}\label{sec3:eq7}
\nonumber &\frac{\partial u_k(t,m)}{\partial t}-\frac{1}{2}\int\limits_{0}^{m} g(m-n,n) [u_0(t,m-n)u_{k-1}(t,n)+u_1(t,m-n)u_k-2(t,n)+\hdots+u_1(t,m-n)\\
&u_1(t,n)] dn+\int\limits_{0}^{\infty}g(m, n)[u_{k-1}(t,m)u_0(t,n)+u_{k-2}(t,m)u_1(t,n)+\hdots
+u_1(t,m)u_1(t,n)]dn=0
\end{align}
The $n$th order approximation to the  solution is defined  as $$\phi_n=\displaystyle \sum_{k=0}^{n} u_k.$$
Note that as ${n\rightarrow \infty}$, we have the limit $\phi=\displaystyle \sum_{k=0}^{\infty} u_k$, which is an analytical solution of  \eqref{sec1:eq1}.
In the next section \ref{s:brkmodel}, the HPM will used to obtain analytical solution of breakage equation for various breakage kernels with different initial conditions.

\section{Application of the HPM in fragmentation equation}
\label{s:brkmodel}
\subsection{Case I: Random binary fragmentation with linear selection function}
In random binary fragmentation the distribution of fragments is given by
\begin{equation}
    k(m|n) = \frac{2}{n} .
\end{equation}
Here we obtain the solution with linear fragmentation rate, $a(m) = m$.  The general form of the factors $u_k(m,t)$ is given by Eq.\ (\ref{pk}), which now becomes
\begin{equation}
     \frac{\partial u_k(t,m)}{\partial t} =
     x u_{k-1}(t,m) - 2\int_0^\infty  \frac{u_{k-1}(t,n)}{n} dn.
\end{equation}
Below we obtain the solution for monodisperse and exponential initial distribution.

\subsubsection{Monodisperse initial conditions}
We start with $u_0=f_0=\delta(m-1)$ and obtain the functions $u_k$ recursively:
\begin{align}
    & u_0 = \delta(m) \\
    & u_1(m,t) = 2 t \theta (1-m)-t x \delta (m-1) \\
    &\vdots \nonumber\\
    & u_k(m,t) = \frac{t^2 \left(1-m\right) (-t m)^{k-2}}{(k-2)!}\theta (1-m),
\end{align}
where $\theta$ is the Heaviside step function. The distribution is given by the infinite series
\begin{equation}
    c(m,t) = \sum_{k=0}^\infty \frac{t^2 \left(1-m\right) (-t m)^{k-2}}{(k-2)!}\theta (1-m),
\end{equation}
As shown in the Appendix, this is equal to
\begin{equation}
\label{austin:sol}
    c(m,t) =  e^{-tm}\left(\delta(m-1)+ \theta(a-m)\left(2t+t^2(1-m)\right) \right)
\end{equation}
This is the same as the result given by Ziff and McGrady \citep{Ziff:JPMG85}.

\subsubsection{Exponential initial condition}
Starting with $u_0(t,m) = f_0(t,m)=e^{-m}$, the $u_k$ are
\begin{align}
   & u_0(m,t) = e^{-m}, \\
   & u_1(m,t) = t \left(-e^{-m}\right) (m-2),\\
   & \vdots \nonumber\\
   & u_k(m,t) =\frac{(-m)^k t^{k-2} \left(-2 k t+(k-1) k+t^2\right) }{k!} e^{-t}   .
\end{align}
The size distribution in this case is given by
\begin{equation}
    c(m,t) = \sum_{k=0}^\infty \frac{(-m)^k t^{k-2}
             \left(-2 k t+(k-1) k+t^2\right) }{k!} e^{-t}
           \to (1+t)^{2}e^{-m(1+t)} .
\end{equation}
This is the same as the solution given by Ziff and McGrady \citep{Ziff:JPMG85}.

\subsection{Case II: Random binary fragmentation with quadratic selection function}
In this model the breakup rate is a quadratic function of size, $a(m)=m^2$ and the fragment size distribution is again given by $g(m|n)=2/n$ and the initial condition is exponential. The general recursion for $u_k(m,t)$ is
\begin{equation}
     \frac{\partial u_k(t,m)}{\partial t} = \int_0^t
     m^2 u_{k-1}(t,m) - 2\int_0^\infty  \frac{u_{k-1}(t,n)}{n^2} dn.
\end{equation}
We obtain solutions for monodisperse and exponential initial conditions:
\subsubsection{Monodisperse initial condition}
With $u_0(m,t)=\delta(m-1)$ the functions $u_k$ are:
\begin{align}
   & u_0(m,t) = \delta(m-1) \\
   & u_1(m,t) = 2 t \theta (1-m)-t x^2 \delta (m-1)\\
   &\vdots\nonumber\\
   & u_k(m,t) = \frac{\delta (t-1) \left(-t^2 x\right)^k}{k!}+\frac{2 m \theta
   (1-t) \left(-t^2 m\right)^{k-1}}{(k-1)!}
\end{align}
The size distribution is
\begin{multline}
    c(m,t) = \sum_{k=0}^\infty
             \left\{ \frac{\delta (t-1) \left(-t^2 x\right)^k}{k!}+\frac{2 x \theta
             (1-t) \left(-t^2 x\right)^{k-1}}{(k-1)!} \right\}
   \\
   \to e^{-tm^2}\left(\delta(m-1)+2at \theta(1-m)\right) .
\end{multline}
This again recovers the solution of Ziff and McGrady \citep{Ziff:JPMG85}

\subsubsection{Exponential initial condition}
Starting with $c_0(m) = e^{-m}$ and solving recursively we obtain
\begin{equation}
    u_k(m,t) = \dfrac{(-t)^k m^{2 k-2}}{k!} \left(-2 k m-2 k+m^2\right) e^{-m} .
\end{equation}
The solution is constructed by computing the infinite series,
\begin{multline}
    c(m,t) =\sum_{k=0}^\infty \dfrac{(-t)^k m^{2 k-2}}{k!} \left(-2 k m-2 k+m^2\right) e^{-m}
    \\
           \to e^{-tm^2}\left(\delta(m-1)+2 t \theta(1-m)\right)  .
\end{multline}
This is the same as the result obtained by Ziff and McGrady \citep{Ziff:JPMG85}

\subsection{Case III: $k$-nary fragmentation with power-law rate}
In this model the selection function is $a(m) = m^\alpha$ and the fragment size distribution is given by
\begin{equation}
    k(m|n) = \frac{\alpha}{n} \left(\frac{m}{n}\right)^{\alpha-2}
\end{equation}
and number of fragments
\begin{equation}
    \bar f = \frac{\alpha}{\alpha -1} .
\end{equation}
By varying $\alpha$ between 1 and 2 the number of fragments ranges from $\infty$ to 2. With $\alpha = 2$, in particular, this model defaults to the random distribution of fragments with quadratic breakage rate. Equation (\ref{pk}) in this model now is
\begin{equation}
   \frac{\partial u_k(m,t)}{\partial t} = x^\alpha u_{k-1}(t,m)
       -  \alpha \int_m^\infty  n m^{\alpha-2} u_{k-1}(t,n) dn.
\end{equation}
Using $u_0(m,t) = e^{-m}$ the functions $u_k$ are
\begin{equation}
    u_k(m,t) = \dfrac{(-t)^k m^{\alpha  k-2}}{k!}
               \left(-\alpha k-\alpha  k m+m^2\right) e^{-m} .
\end{equation}
Inserting this result into the infinite series in Eq.\ (\ref{series}) we obtain
\begin{multline}
   c(t,m) = \sum_{k=0}^\infty
            \dfrac{(-t)^k m^{\alpha  k-2}}{k!}
               \left(-\alpha  k-\alpha  k m+m^2\right) e^{-m} \\
         \to
    e^{-m-t m^\alpha} (1+ \alpha t (m^{\alpha-2}+m^{\alpha-1})) ,
\end{multline}
which gives the size distribution at all times. Now, HPM will be apply to solve various models of aggregation equation.


\section{Application of the HPM in aggregation equation}
\label{s:aggmodel}

\subsubsection{Constant aggregation kernel $g(m, n)=1$}
The distribution of aggregates for constant aggregation kernel is
\begin{align*}
g(m, n)=1.
\end{align*}
The general form of the factors $u_k(m, t)$ is given by Eq. \eqref{sec3:eq7}, which now becomes
\begin{align}
\nonumber &\frac{\partial u_k(t,m)}{\partial t}-\frac{1}{2}\int\limits_{0}^{m} g(m-n, n) [u_0(t,m-n)u_{k-1}(t,n)+u_1(t,m-n)u_k-2(t,n)+\hdots+u_1(t,m-n)\\&
\nonumber u_1(t,n)] dn+\int\limits_{0}^{\infty}g(m, n)[u_{k-1}(t,m)u_0(t,n)+u_{k-2}(t,m)u_1(t,n)+\hdots+u_1(t,m)u_1(t,n)]dn=0
,\\&\hspace{12cm}k=1,2,3...
\end{align}
The solution will be obtained for exponential initial condition, i.e., $u_0(m)=e^{-m}$. So, let us initialize the solution with $u_0 =e^{-m}$ and obtain the functions $u_k$ recursively
\begin{align*}
u_0&=e^{-m},\\
u_1&= \frac{t}{2^1 1!}  (m-2)e^{-m},\\
u_2&=\frac{t^2}{ 2^2 2!}  \left(m^2-6 m+6\right) e^{-m} ,\\
\vdots\\
u_{k}&= e^{-m} \frac{(k+1)!\; t^k}{ 2^{k} } \sum_{r=0}^{k}\frac{(-1)^{r}\; m^{k-r}}{r!(k-r+1)!(k-r)!}
\end{align*}
The partial sum of the series solution is defined   obtained as
\begin{align*}
 u(t,m)\approx \phi_n(t,m)=\sum_{k=0}^{n} u_k.
 \end{align*}
 As $n\longrightarrow \infty$, we obtain
 \begin{align*}
  u(t,m)=e^{-m} \sum _{k=0}^{\infty } \frac{t^k }{2^{k} k!} HypergeometricU[-k,2,m]
 \end{align*}
The  exact solution is
$u(t,m)=n^{2}(t)e^{-n(t)m},$ where $n(t)=\frac{2}{t+2}$ as  in \cite{ranjbar2010numerical}.

\subsection{Sum aggregation kernel $g(m, n)=m+n$}
The distribution of aggregates for sum aggregation kernel is defined as
\begin{align*}
g(m, n)=m+n
\end{align*}
The analytical  solution can be found in \cite{kumar2006numerical} which is $$u(t,m)= \frac{(1-\tau ) e^{-(1+\tau)m}}{m\sqrt{\tau}}I[2m\sqrt{\tau}],$$ where
$\tau =1-e^{-t}$ and  $I[m]$ is the modified Bessel function of the first kind.
The general form of the factors $u_k(m, t)$ is given by Eq. \eqref{sec3:eq7}, which now becomes
\begin{align}
\frac{\partial u_k}{\partial t}=\displaystyle  \frac{1}{2}\int\limits_{0}^{m} m \;  H^{1}_{k-1} dn-\int\limits_{0}^{\infty} (m +n) \; H^{2}_{k-1}dy,~~~~k=1,2,3...
\end{align}

\subsubsection{Exponential initial condition $u_0(m)=e^{-m}$}
\begin{figure}
\begin{center}
\includegraphics[width=4in]{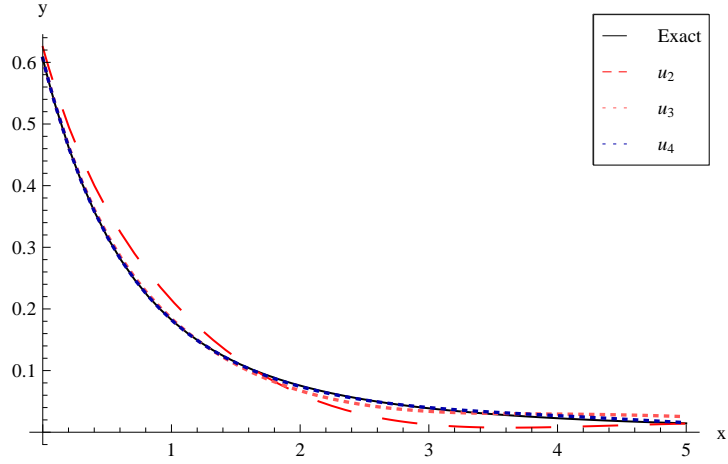}
\end{center}
\caption{Comparison of Exact solution with approximate solution.}
\label{sumkernel}
\end{figure}

We initialize with $u_0 = e^{-m}$ and obtain the functions $u_k$ recursively
\begin{align*}
u_0&= e^{-m},\\
u_1&=\frac{1}{2} e^{-m} t \left(m^2-2 m-2\right),\\
u_2&=\frac{1}{12} e^{-m} t^2 \left(6+18 m-3 m^2-6 m^3+m^4\right),\\
u_3&=\frac{1}{144} e^{-m} t^3 \left(-24-168 m-60 m^2+120 m^3+12 m^4-12 m^5+m^6\right)\\
u_4&=\frac{1}{2880} e^{-m} t^3 \left(120+1800 m+2100 m^2-1800 m^3-1180 m^4+360 m^5+70 m^6-20 m^7+ m^8\right)\\
\vdots\\
\end{align*}
It is difficult to find the generalized term for the solution of sum kernel. So, we approximate the solution by taking the summation of some terms. The comparison of exact solution with the approximate solution is represented in Figure \ref{sumkernel}. The comparison reveals that by considering the sum of four terms in the series, exact solution can be approached.

\subsection{Product aggregation kernel $g(m, n)=mn$}
The distribution of aggregates for product aggregation kernel is
\begin{align*}
g(m, n)=mn
\end{align*}
The general form of the factors $u_k(m, t)$ is given by Eq. \eqref{sec3:eq7}, which now becomes
\begin{align}
\nonumber &\frac{\partial u_k(t,m)}{\partial t}-\frac{1}{2}\int\limits_{0}^{m} (m-n)n [u_0(t,m-n)u_{k-1}(t,n)+u_1(t,m-n)u_k-2(t,n)+\hdots+u_1(t,m-n)\\&
\nonumber u_1(t,n)] dn+\int\limits_{0}^{\infty}mn[u_{k-1}(t,m)u_0(t,n)+u_{k-2}(t,m)u_1(t,n)+\hdots+u_1(t,m)u_1(t,n)]dn=0
,\\&\hspace{12cm}k=1,2,3...
\end{align}
The solution will be obtained for two different initial conditions, i.e., for $u_0(m)= e^{-m}$ and $\displaystyle{u_0(m)= \frac{e^{-m}}{m}}$.

\subsubsection{Exponential initial condition, $u_0(m)= e^{-m}$}
Its analytic solution follows \cite{kumar2006numerical} as  $$u(t,m)=\sum _{k=0}^{\infty } e^{-(t+1)m} \frac{t^k m^{3k}}{(k+1)! \; \Gamma(2k+2)}$$

We start with $u_0 =e^{-m}$ and obtain the functions $u_k$ recursively
\begin{align*}
u_0&= e^{-m},\\
u_1&=\frac{(tm)^{1} }{2! 3!}  (-12+m^2)e^{-m},\\
u_2&=\frac{(tm)^{2}}{ 3!5! }  \left(360-60 m^2+m^4\right)e^{-m},\\
\vdots\\
u_{k}&=  2e^{-m}(mt)^k \sum_{r=0}^{k}\frac{(-1)^{r}\; m^{2k-2r}}{r!(k-r)!\; (2k-2r+2)!}
\end{align*}
The partial sum of the series solution is defined   obtained as
\begin{align*}
 u(t,m)\approx \phi_n(t,m)=\sum_{k=0}^{n} u_k.
 \end{align*}
Now, the solution for another initial condition.
\subsubsection{When $u_0(m)= \frac{e^{-m}}{m}$}
The analytical solution can be found in \cite{ranjbar2010numerical}
\begin{align*}
 u(t,m)=\frac{I[2m\sqrt{t}]}{m^{2}\sqrt{t}}\; e^{-Tm},
\end{align*}
where
$T = \left\{
  \begin{array}{ll}
   1+t , & \quad \text{ $t\leq1$ },\\
    2\sqrt{t}, &  \quad \text{otherwise}.
  \end{array} \right.$\\

Here $I$ is the modified Bessel's first kind function:
\begin{align}
I[m]=\frac{1}{\pi}\int\limits_{0}^{\pi}e^{(m \cos\theta)}\cos\theta d\theta.
\end{align}
We initialize with $u_0 = \frac{e^{-m}}{m}$ and obtain the functions $u_k$ recursively
\begin{align*}
u_0&=\displaystyle \frac{e^{-m}m^{-1}}{1},\\
u_1&=\displaystyle\frac{t}{2 (1!)^2}  (m-2)e^{-m} ,\\
u_2&=\displaystyle\frac{t^2 m}{3 (2!)^2}  \left(m^2-6 m+6\right)e^{-m} ,\\
\vdots\\
u_{k}&= e^{-m} x^{k-1 }t^k \sum_{r=0}^{k}  \frac{(-1)^{r} m^{k-r}}{r!(k-r+1)!(k-r)!}.
\end{align*}
The partial sum of the series solution is obtained as
\begin{align*}
 u(t,m)\approx \phi_n(t,m)=\sum_{k=0}^{n} u_k.
 \end{align*}
 As $n\longrightarrow \infty$, we obtain
\begin{align*}
e^{-m} \sum _{k=0}^{\infty } \frac{t^k m^{k-1} }{k! (k+1)!} HypergeometricU[-k,2,m]
 \end{align*}

\section{Discussion}\label{s:discussion}
Cases 1 and 2, which describe random binary breakup with linear and quadratic rate functions essentially form the ``classical'' model in fragmentation. Their solutions have been obtained previously in the literature multiple times by various methods, including probabilistic treatment \cite{Ziff:JPMG85}, and solution of the equivalent discrete problem \citep{Montroll:JCP40}. The kernel of Case 3 is a special form of a kernel whose general form is
\begin{equation}\label{austin}
    k(m|n) =
    (1-\psi ) \lambda  \frac{m^{\lambda -1}}{n^{\lambda }} +
   \psi   \gamma  \frac{m^{\gamma -1}}{n^{\gamma }} .
\end{equation}
This kernel was introduced by Austin \citep{Austin:PT72,Klimpel:IJMP77} to model grinding. In cumulative form it is the sum of two power-law distributions, of which one is the dominant distribution for large $m$ (Austin gave this kernel with different coefficients but the two forms can be shown to be the same under appropriate relationships among the coefficients). The same kernel has been the subject of theoretical investigations of the scaling behavior of the fragmentation equation \citep{Ziff:JPMG91} and of the so-called shattering transition \citep{Ziff:M86,Ernst:JPMG93}, a behavior analogous but opposite to gelation that is manifested by the presence of a finite fraction of the mass of the system in particles of zero size.  Ziff and McGrady \citep{Ziff:JPMG85,McGrady:PRL87} gave the solution of the Austin kernel with power-law breakage rate for monodisperse initial conditions. Oukouomi Noutchie and Doungmo Goufo obtained the same solution by Laplace transform \citep{Oukouomi-Noutchie:MPE14}.
\begin{figure}
\begin{center}
\includegraphics[width=3.25in]{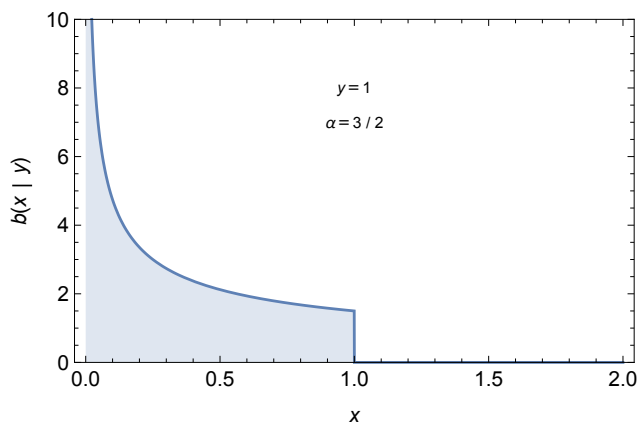}
\end{center}
\caption{The distribution $k(m|n)$ for $\alpha=3/2$, $n=1$. }
\label{fig_bxy}
\end{figure}
In principle, given the linearity of the fragmentation equation the solution for any other initial condition can be reconstructed from the solution to the monodisperse problem. The monodisperse solution, however, employs the confluent hypergeometric function and its manipulation is not trivial. The solution obtained here in Eq.\ (\ref{austin:sol}) for exponential conditions is compact and represents a new result. This solution obeys the scaling derived by Ziff and McGrady \citep{Ziff:JPMG85}, as can be easily demonstrated. Dividing Eq.\ (\ref{austin:sol}) by $m^2/\alpha$ and setting $m=(z/t)^{1/\alpha}$ we find
\begin{equation*}
\frac{\alpha c(m,t)}{m^2} =
    \frac{1}{\alpha}\left\{
    t^{-2/\alpha } e^{t^{-1/\alpha }
       \left(
       - z^{1/\alpha }\right)-z} \left(\alpha  t^{1/\alpha } z^{(\alpha+1)/\alpha}
       + \alpha  z t^{2/\alpha }
       + z^{2/\alpha }
       \right)
       \right\} ,
\end{equation*}
whose limit at long $t$ is
\begin{equation*}
           \frac{\alpha c(m,t)}{m^2} \to z e^{-z} \equiv \Phi(z) ,
\end{equation*}
and $\Phi(z)$ is the scaling function. The scaling limit of the distribution is
\begin{equation}
    c(m,t) \sim \alpha m^2\Phi(t m^a) .
\end{equation}
%
\begin{figure}
\begin{center}
\includegraphics[width=3.25in]{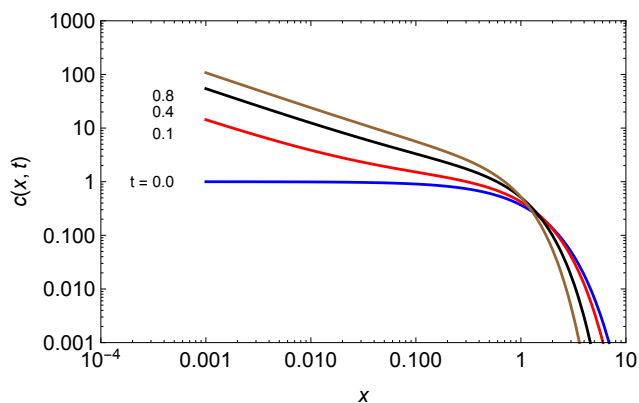}
\end{center}
\caption{Evolution of the size distribution for $\alpha=3/2$. }
\label{fig_dstr}
\end{figure}
%
The scaling function given by Ziff and McGrady is $\Phi(z) = z^{\gamma/\alpha} e^{-z}\Gamma(\gamma/\alpha)$ (the Filippov model in Table 1 of Ref.\ \citep{Ziff:JPMG85} ), which for $\gamma=\alpha$ reduces to the result obtained here.

To illustrate these results we examine the solution for $\alpha = 4/3$. The fragmentation rate is $a(m)=m^{4/3}$ and the distribution of fragments is
\begin{equation}
    k(m|n) = \frac{4}{3 m^{2/3} n^{3/2}},
\end{equation}
with an average number of fragments $f=4$. This distribution is very asymmetric with a large fraction of the population in the small size region (Fig.\ \ref{fig_bxy}). There is no shattering, however, because the fragmentation rate decreases as the size approaches zero (indeed our solution does include shattering as the condition $f\geq 2$ constrains $\alpha$ to be positive). The size distribution is shown in Fig.\ \ref{fig_dstr}. After an initial transition the distribution in log-log coordinates reaches an invariant form that is simply translated to lower average sizes. This is a manifestation of the scaling limit.

\begin{figure}
\begin{center}
\includegraphics[width=3.25in]{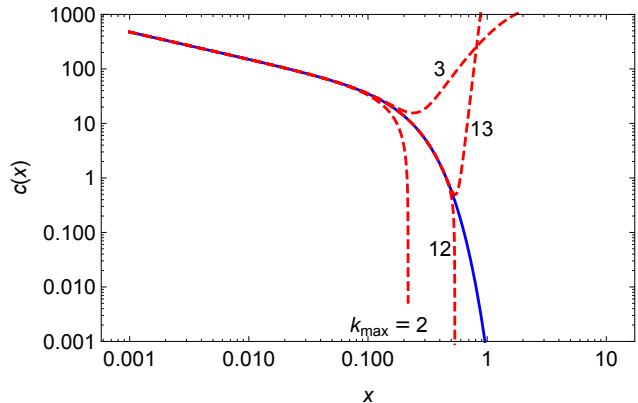}
\end{center}
\caption{Truncated series to $k_\text{max}$ for $\alpha=3/2$, $t=10$. The dashed lines are truncated approximations of the full solution (sold line). }
\label{fig_series}
\end{figure}

Finally, to examine the convergence of the series we calculate the solution by retaining only a finite number of terms $k$. Figure \ref{fig_series} compares truncated approximations of the full solutions for $\alpha=3/2$ at $t=10$ and $k_\text{max}=2,3,12$ and $13$ ($k_\text{max}$ is the maximum order of the term retained in the series). The general feature of the truncation is that the approximation is excellent for $0<m<m_\text{max}$, where $m_\text{max}$ is the size beyond which the approximation breaks down and is pushed to larger sizes as the order of the truncation increases. Even-order truncations provide stable approximations approximations of the size distribution and it is interesting to point out that in this particular case, even a very low order approximation such as $k_\text{max}=2$ represents the distribution very well over two decades in $c(m)$.

\section{Conclusions}
The homotopy perturbation method has been successfully applied to the fragmentation as well as aggregation population balance equation. The HPM provides a systematic methodology to solve for the size distribution and was shown to reproduce all known solutions and to lead to new results for a special case of the Austin kernel and also reproduce analytical solution for various models of fragmentation and aggregation.


\bibliography{hpm_paper}
\bibliographystyle{apsrev}
\end{document}